\documentclass[12pt]{article}
\usepackage{amssymb}
\usepackage{amsmath}
\usepackage{endnotes}
\usepackage[doublespacing]{setspace}

\input {tcilatex}
\begin{document}

\title{Nonparametric functionals as generalized functions }
\author{Victoria Zinde-Walsh \thanks{%
The support of the \textit{Fonds\ qu\'{e}becois de la recherche sur la soci%
\'{e}t\'{e} et la culture} (FRQSC) is gratefully acknowledged. } \\
\\
McGill University and CIREQ\\
victoria.zinde-walsh@mcgill.ca\\
(514) 398 4834}
\maketitle
\date{}

\begin{center}
\bigskip \pagebreak
\end{center}

Running head: Nonparametric functionals

Victoria Zinde-Walsh

Department of Economics, McGill University

855 Sherbrooke Street West,

Montreal, Quebec, Canada

H3A 2T7

\begin{center}
{\LARGE Abstract}
\end{center}

The paper considers probability distribution, density, conditional
distribution and density and conditional moments as well as their kernel
estimators in spaces of generalized functions. This approach does not
require restrictions on classes of distributions common in nonparametric
estimation. Density in usual function spaces is not well-posed; this paper
establishes existence and well-posedness of the generalized density
function. It also demonstrates root-n convergence of the kernel density
estimator in the space of generalized functions. It is shown that the usual
kernel estimator of the conditional distribution converges at a parametric
rate as a random process in the space of generalized functions to a limit
Gaussian process regardless of pointwise existence of the conditional
distribution. Conditional moments such as conditional mean are also be
characterized via generalized functions. Convergence of the kernel
estimators to the limit Gaussian process is shown to hold as long as the
appropriate moments exist.

\section{Introduction}

A probability distribution function, $F,$ that corresponds to a Borel
measure on a Euclidean space $R^{k}$ (or its subspace) is always defined in
the space of bounded functions. It can be viewed as the right-hand side of
an integral equation: 
\begin{equation}
I(f)=F;  \label{inverseprobl}
\end{equation}%
where the density represents the solution to the inverse problem 
\begin{equation}
f=\partial ^{k}F.  \label{solution}
\end{equation}%
Here $I$ represents an integration operator for $R^{k}$: $I\left( f\right)
(x)=\int_{-\infty }^{x_{1}}...\int_{-\infty }^{x_{k}}f\left( w\right)
dw_{1}...dw_{k}$ and $\partial ^{k}=\frac{\partial ^{k}}{\partial
x_{1}...\partial x_{k}}$ the inverse differentiation operator.

When does the solution to the inverse problem exist?

In the usual approach the integral operator $I$ is assumed to operate on the
space of integrable functions, e.g. $L_{1}$ (absolutely integrable
functions) or $L_{2}$ (square integrable functions), - see e.g. Devroye and
Gyorfi (1985), Carrasco, Florens, Renault (2007). The operator $I$ maps
density functions in $L_{1}$ into the space of absolutely continuous
distribution functions. In this case the inverse operator $\partial ^{k}$ is
defined and the inverse problem has a unique solution.

The property of well-posedness requires that the solution continuously
depend on the right-hand side function, in other words, if distribution
functions are close, the corresponding densities should be close as well.
However, in spaces of integrable functions the inverse problem is not
well-posed: while the operator $I$ is continuous on $L_{1}$ (or another $%
L_{p}$ space) the inverse operator $\partial ^{k}$ is not. The example below
(from Zinde-Walsh, 2011) illustrates lack of well-posedness.

Example. \textit{Consider the space }$D\left( \left[ 0,1\right] \right) $ 
\textit{of univariate absolutely continuous distribution functions on the
interval }$\left[ 0,1\right] $\textit{\ in the uniform metric: the distance
between two distributions, }$F_{1},$ \thinspace $F_{2}$ is $d(F_{1},F_{2})=%
\underset{x\in \lbrack 0,1]}{\max }\left\vert F_{1}(x)-F_{2}(x)\right\vert ;$
\textit{this is the image space of the operator }$I(\cdot )$\textit{\
defined on }$L_{1}\left( \left[ 0,1\right] \right) .$

\textit{Denote by }$\left[ v\right] $\textit{\ the integer part of }$v,$%
\textit{\ that is the largest integer that is }$\leq v.$\textit{\ Let }$%
I\left( x\in A\right) $\textit{\ denote the indicator function of set }$A,$%
\textit{\ that equals }$1$\textit{\ if }$x$\textit{\ is in }$A,$\textit{\
zero otherwise. With }$\varepsilon =\frac{\bar{\varepsilon}}{2}$\textit{\
define densities}%
\begin{eqnarray*}
f_{1}(x) &=&2\dsum\limits_{m=0}^{\left[ \frac{\varepsilon ^{-1}+1}{2}\right]
-1}I\left( x\in \lbrack 2m\varepsilon ,\left( 2m+1\right) \varepsilon
\right) ); \\
f_{2}\left( x\right) &=&2\dsum\limits_{m=0}^{\left[ \frac{\varepsilon ^{-1}+1%
}{2}\right] -1}I\left( x\in \lbrack (2m+1)\varepsilon ,\left( 2m+2\right)
\varepsilon \right) ).
\end{eqnarray*}

\textit{The densities }$f_{1}$\textit{\ and }$f_{2\text{ }}$\textit{have
supports that do not intersect, it is easily seen that at each point they
differ by 2: }$\left\vert f_{1}\left( x\right) -f_{2}\left( x\right)
\right\vert =2$\textit{; it follows that the }$L_{1}\left( \left[ 0,1\right]
\right) $ \textit{difference between them is }$2.$ \textit{The corresponding
distributions are }$F_{1}=I(f_{1})$\textit{\ and }$F_{2}=I(f_{2}).$\textit{\
It is easy to establish by integration that }%
\begin{equation*}
\underset{x\in \lbrack 0,1]}{\max }|F_{1}\left( x\right) -F_{2}\left(
x\right) |\leq 2\varepsilon =\bar{\varepsilon},
\end{equation*}%
\textit{and thus the inverse operator is not continuous.}

Thus although a solution to the inverse problem in the $L_{1}$ space exists
for absolutely continuous distributions, the problem is not well-posed. 

By contrast, in the appropriate space of generalized functions the solution
to the density problem exists without any restrictions on the distribution
function and is well-posed; as proved in section 2 below this follows from
the known properties of generalized functions. The fact that generalized
functions can be useful when non-differentiability prevents the use of
Taylor expansions was discussed e.g. in Phillips (1991) for LAD estimation,
and continued in some econometric literature that followed.

The statistical inverse problem is solved often with a kernel density
estimator. Consider a random sample of observations from a distribution $F,$ 
$\{x_{i}\}_{i=1}^{n},$ $x_{i}\in R^{k}.$ With a chosen kernel function, $K$
and bandwidth (vector) $h$ the estimator is 
\begin{equation}
\widehat{f(x)}=\frac{1}{n\Pi _{j=1}^{k}h_{j}}\sum_{i=1}^{n}K(\frac{x_{i}-x}{h%
}),  \label{kernel}
\end{equation}%
where $h$ has components $h_{1},...h_{k}$ and $K(\frac{x_{i}-x}{h})$ is a
multivariate function with the argument $\left( \frac{x_{i}-x}{h}\right)
=\left( \left( \frac{x_{i1}-x_{1}}{h_{1}}\right) ,...,\left( \frac{%
x_{ik}-x_{k}}{h_{k}}\right) \right) $. We shall proceed with the following
assumption on the kernel.

\textbf{Assumption 1 (kernel).}

\textit{(a). }$K(w)$\textit{\ is an ordinary bounded function on }$R^{k};$%
\textit{\ }$\int K(w)dw=1;$

\textit{(b). Support of }$K$\textit{\ belongs to }$[-1,1]^{k};$

\textit{(c). }$K(w)$\textit{\ is an }$l-$\textit{th order kernel: for }$%
w=(w_{1},...w_{k})$\textit{\ the integral }%
\begin{equation*}
\int w_{1}^{j_{1}}...w_{k}^{j_{k}}K(w)dw_{1}...dw_{k}\left\{ 
\begin{array}{cc}
=0 & \text{\textit{if }}\mathit{j}_{1}+...+\mathit{j}_{k}\mathit{<l;} \\ 
<\infty & \text{\textit{if }}\mathit{j}_{1}+...+\mathit{j}_{k}\mathit{=l.}%
\end{array}%
\right.
\end{equation*}

The finite support and boundedness assumptions can be relaxed and are
introduced to simplify assumptions and derivations; $K$ is not restricted to
be symmetric or non-negative.

Denote by $\bar{K}$ the integral of the kernel function, then 
\begin{equation}
\widehat{F(x)}=\frac{1}{n}\sum_{i=1}^{n}\bar{K}(\frac{x_{i}-x}{h})
\label{distrib}
\end{equation}%
is an estimator of the distribution function, $F(x).$ The properties of
these estimators depend on $K$ and $h$ and are well established (Azzalini,
1981). Generally for $h\rightarrow 0$ as $n\rightarrow \infty $ with $%
nh\rightarrow \infty $, $\widehat{F(x)}$ is a root-n consistent and
asymptotically Gaussian estimator of $F(x)$ at any point of continuity; the
uniform norm of the difference, $\sup \left\vert \widehat{F(x)}%
-F(x)\right\vert ,$ converges to zero.

Known convergence properties of $\ \widehat{f(x)}$ are more complicated;
they rely on assumptions about the existence and smoothness of the density, $%
f(x);$ the convergence rate is slower than root-n and depends on the order
of the kernel and the rate of the bandwidth $h\rightarrow 0$ (Pagan and
Ullah, 1999). As shown in Examples 3-5 \ in Zinde-Walsh (2008), the
estimator $\widehat{f(x)}$ fails to converge pointwise if the distribution
is not absolutely continuous (e.g. at a mass point or for a fractal
measure); of course, in those cases density itself cannot be defined
pointwise and exists only as the solution, $f$ in $\left( \ref{solution}%
\right) $ to the inverse problem in the space of generalized functions.

When considered in the space of generalized functions the estimators, $%
\widehat{f},$ are viewed as random continuous linear functionals on spaces
of well-behaved functions where convergence to generalized derivatives of
distribution functions (solutions to the inverse problem) can be established
without any assumptions on the underlying distribution. Moreover,
convergence of kernel estimators can be faster and even at parametric rates.
This result has features common to other results on convergence of random
functionals of density as discussed, e.g. in Anderson et al (2012) and is
derived here in section 3. This result relies on the derivation of the rate
of bias in generalized functions that was provided in Zinde-Walsh (2008) but
gives the derivation of the covariance functional that corrects the one in
that paper.

Conditioning is somewhat awkward and there are many different ways to
streamline the representation of conditional measures and distribution
functions (Chiang and Pollard, 1997, Pfanzagl, 1979 among others). Here we
focus on the distribution function $F\left( x,y\right) $ on $R^{d_{x}}\times
R^{d_{y}}$ and distribution of $y\in R^{d_{y}}$ conditional on $x\in
R^{d_{x}}.$ In this case typically the conditional distribution $F_{y|x}$
function is represented via a fraction $\frac{\partial ^{d_{x}}F(x,y)}{%
f_{x}\left( x\right) },$ where the differentiation operator is applied to
the $x$ argument of $F\left( x,y\right) $ and $f_{x}\left( x\right) $
represents the density of the marginal distribution. Of course such a
representation makes stringent requirements on the smoothness of the
appropriate functions. Here the case of an arbitrary continuous conditioning
distribution is considered without requiring differentiability; it is shown
that for this case the conditional distribution and conditional density have
a straightforward representation as generalized functions on appropriate
spaces. The representation is in terms of functionals involving the
conditioning distribution (rather than the conditioning variable) as an
argument; this representation avoids the nonlinearity introduced by the
denominator. When the usual representation holds, a simple correspondence
between the two representations is established. Conditional density, $f_{y|x}
$ is defiend as a generalized derivative of the conditional distribution
generalized function.

The convergence of the usual kernel estimator of the conditional
distribution is known under smoothness assumptions (Pagan and Ullah, 1999,
Li and Racine, 2007) and utilizes the properties of the kernel density
estimator; the density appears in the denominator of the statistic requiring
some support assumptions and possibly regularization to converge. Here the
root-n convergence of the kernel estimator to a limit Gaussian process in
generalized function space is established without any extra restrictions on
the distribution.

An interpretation of a conditional moment function is provided here in the
space of generalized functions, thus again without any restriction beyond
continuity of conditioning distribution. For estimators, such as for
conditional mean kernel estimator the asymptotic properties are established,
the result is then that root-n convergence in generalized functions obtains
for the kernel estimator without any restrictions on smoothness of
distribution functions.

The theoretical results of this paper extend the usual representation of the
density, conditional distribution and density and conditional moments to
situations where these may not exist in an ordinary sense. The advantage
that this approach provides is its generality. On the other hand, the
topology in the spaces of generalized functions is weak and well-posedness
does not imply convergence in norm.

The asymptotic results provide a general approach, so that when the usual
assumptions may fail there is still a sense in which consistency holds.
Moreover a root-n convergence rate obtains, again as a consequence of the
weak topology with no guarantee of good convergence in norm. The practical
advantage is in the possibility of utilizing the generalized random process
and its limit process for inference without making any restrictions on the
distribution.

\section{Density as solution to a well-posed inverse problem in the space of
generalized functions}

For the definitions and results pertaining to spaces of generalized
functions the main references are to books by Schwartz (1966) Gel'fand and
Shilov (1964). A useful summary is in Zinde-Walsh (2008, 2012); the main
definitions follow.

Consider a space of well-behaved "test" functions, $D_{\infty }\left(
R^{k}\right) $ of infinitely differentiable functions with bounded support,
or any of the spaces $D_{m}\left( R^{k}\right) $ of $m$ times continuously
differentiable functions (with bounded support); sometimes the domain of
definition can be an open subset $W$ of $R^{k},$ typically here $W=\left(
0,1\right) ^{k}.$ Denote the generic space by $D\left( W\right) ;$
convergence in $D\left( W\right) $ is defined as follows: a sequence $\psi
_{n}\in D\left( W\right) $ converges to zero if all $\psi _{n}$ are defined
on a common bounded support in $W$ and $\psi _{n\text{ }}$ as well as all
the $l-th$ order derivatives (with $l\leq m$ for $D_{m}$ or all $l<\infty $
for $D_{\infty })$ converge pointwise to zero. The space of generalized
functions is the dual space, $D^{\ast },$ the space of linear continuous
functionals on $D\left( W\right) $ with the weak topology: a sequence of
elements of $D^{\ast }$ converges if the sequence of values of the
functionals converges for any test function from $D\left( W\right) .$ The
usual notation is to write the value of the functional $f$ applied to a test
function $\psi \in D\left( W\right) $ as $\left( f,\psi \right) ;$ then a
sequence $f_{n}$ converges to $f$ if for any $\psi $ convergence $\left(
f_{n},\psi \right) \rightarrow \left( f,\psi \right) $ holds.

Assume that functions in $D\left( W\right) ;$ $W\subseteq R^{k}$ are
suitably differentiable$,$ e.g. at least $k$ times continuously
differentiable. Then for any $\psi \in D\left( W\right) ,$ and $F\in D^{\ast
}$ define a generalized derivative $f\in D^{\ast };$ $f=\frac{\partial ^{k}}{%
\partial x_{1}...\partial x_{k}}F$ as the functional with values given by$:$ 
\begin{equation}
(f,\psi )=(-1)^{k}(F,\frac{\partial ^{k}\psi }{\partial x_{1}...\partial
x_{k}}).  \label{multidens}
\end{equation}%
If the right-hand side is expressed via a regular locally summable function
as is the case when $F$ is a probability distribution function, then it can
be computed by integration:%
\begin{equation*}
\left( F,\frac{\partial ^{k}\psi }{\partial x_{1}...\partial x_{k}}\right)
=\int ...\int F(x_{1},...,x_{k})\frac{\partial ^{k}\psi (x_{1},...x_{k})}{%
\partial x_{1}...\partial x_{k}}dx_{1}...dx_{k}.
\end{equation*}%
For the function $F$ $\left( \ref{multidens}\right) $\ the functional on the
right-hand side defines the generalized derivative$:$ $f=$ $\frac{\partial
^{k}F}{\partial x_{1}...\partial x_{k}}.$

First consider density as a generalized function on the space $D_{\infty
}\left( W\right) .$

\textbf{Theorem 1.} \textit{The inverse problem }$\left( \ref{inverseprobl}%
\right) $\textit{\ for any cumulative probability distribution function }$F$%
\textit{\ has the solution }$f$\textit{\ defined by }$\left( \ref{multidens}%
\right) $\textit{\ in the space of generalized functions }$D^{\ast }$\textit{%
\ for }$D_{\infty }\left( W\right) $\textit{. The problem is well-posed.
When density exists as an integrable function, }$f(x),$\textit{\ it provides
the generalized function }$f$ \textit{via the value of the corresponding
functional:}%
\begin{equation}
\left( f,\psi \right) =\int ...\int f(x_{1},...,x_{k})\psi
(x_{1},...,x_{k})dx_{1}...dx_{k}.  \label{regdens}
\end{equation}

Proof.

Any distribution function $F$ on $R^{k}$ is a monotone bounded function and
as such is locally integrable on any bounded set; a function like that
represents a regular element in the space of generalized functions, $D^{\ast
},$ for $D_{\infty }\left( W\right) $ defined above. Then $\left( \ref%
{multidens}\right) $ defines $f$ as the generalized derivative of $F,$ the
generalized density function.

The differentiation operator $\partial ^{k}=\frac{\partial ^{k}}{\partial
x_{1}...\partial x_{k}}$on the space of generalized functions $D^{\ast }$ is
defined for any regular function and is a continuous operator (Schwartz,
p.80). Thus the solution $f$ continuously depends on $F$ in these spaces
providing well-posedness.

If density $f$ exists as a regular integrable function, its integral
coincides with the function $F$ and integration by parts of $\left( \ref%
{regdens}\right) $ provides $\left( \ref{multidens}\right) .$ Thus $f$, the
solution to the inverse problem in the space $D^{\ast }$ is consistent with
the solution when it exists as an ordinary function.$\blacksquare $

\textbf{Corollary.}\textit{\ The result of the Theorem applies in the space
of generalized functions on }$D_{m}\left( W\right) ,$ $m\geq k.$

Proof.

Indeed, consider the space $D_{\infty }\left( W\right) \subset D_{k}\left(
W\right) .$ By the theorem the inverse problem provides the density function 
$f$ defined as a linear continuous functional on $D_{\infty }\left( W\right) 
$ via $\left( \ref{multidens}\right) .$ We can extend the functional $f$ to $%
D_{k}\left( W\right) $ as a linear continuous functional. First note that
since $F$ is a regular locally integrable function it represents an element
in $D_{k}^{\ast };$ then define the functional in $D_{k}^{\ast }$ by $\left( %
\ref{multidens}\right) $ for any $\psi \in G_{k},$ denote it $\tilde{f}$ to
distinguish from $f$ defined on $D_{\infty }\left( W\right) .$ This $\tilde{f%
}$ represents a linear continuous functional, so an element in $D_{k}^{\ast
}.$ There is an injective mapping of linear topological spaces $D_{k}^{\ast
}\rightarrow D_{\infty }^{\ast }$ (Sobolev, 1992\ ; in notation there $%
C^{(k)\#}\rightarrow C^{(\infty )\#}),$ thus by this mapping $\tilde{f}$
maps into $f$ and the inverse problem is solved in $D_{k}^{\ast }$ and is
well-posed there .$\blacksquare $

\section{Gaussian limit process for the kernel density estimator in the
space of generalized functions}

We now describe the limit process for the kernel estimator $\left( \ref%
{kernel}\right) $ as $\bar{h}=\underset{1\leq j\leq k}{\max }%
h_{j}\rightarrow 0$ with $n\rightarrow \infty ,$ as a generalized random
process. Such a description was in Zinde-Walsh (2008), but there was an
error in the variance computation that is corrected here. The main result
here is that in the generalized functions space convergence of the kernel
density estimator can be at a parametric rate for a suitable selection of
the kernel and bandwidth; unlike the usual case in the literature this
selection alone provides the result independently of any properties
(smoothness) of the distribution.

Recall that convergence of generalized random functions is defined (see,
e.g. Gel'fand and Vilenkin, 1964 or summary in Zinde-Walsh, 2008) as weak
convergence of random linear continuous functionals on the space $D_{\cdot }$
(for any of the $D_{k},D_{\infty },$ etc. spaces here) that are indexed by
the functions in $D_{\cdot }$: stochastic convergence of random functionals, 
$\hat{f},$ follows from stochastic convergence of random vectors of values
of the functional $\left( \left( \hat{f},\psi _{1}\right) ,...,\left( \hat{f}%
,\psi _{m}\right) \right) ^{\prime }$ for any finite set $\left( \psi
_{1},...,\psi _{m}\right) $ with $\psi _{l}\in D_{\cdot }.$ Thus we need to
consider the behavior of such random vectors.

Theorem 2 in Zinde-Walsh (2008) gives the convergence rate $O(\bar{h}^{l})$
for the generalized bias function of the kernel estimator based on a random
sample and the expression for the bias for $\psi \in D_{l+k}$ and kernel $K$
of order $l:$ 
\begin{equation*}
E\widehat{f}-f\approx O(\overline{h}^{l}),
\end{equation*}%
more specifically for any $\psi $\ the bias functional provides $\left( E%
\widehat{f},\psi \right) -\left( f,\psi \right) =$\ \ \ \ \ \ \ \ \ \ \ \ \
\ 
\begin{eqnarray}
&&(-1)^{l}\dsum\limits_{\Sigma m_{i}=l}\int \dprod\limits_{i=1}^{k}\frac{%
h_{i}^{m_{i}}}{m_{i}!}F(\tilde{x})\frac{\partial ^{l+k}\psi }{\partial
x_{1}^{m_{1i}+1}...\partial x_{k}^{m_{ki}+1}}(\tilde{x})d\tilde{x}\int
K(w)w_{1}^{m_{1}}...w_{k}^{m_{ki}}dw  \label{biasfl} \\
&&+R(h),  \notag
\end{eqnarray}%
where $R(h)=o(\overline{h}^{l});$ if $\psi \in D_{l+k+1}$ then $R(h)=O(%
\overline{h}^{l+1}).$ Note that $\left( f,\psi \right) =E(\psi )$ where
expectation is with respect to the measure given by $F.$

Denote the expression 
\begin{equation*}
(-1)^{l}\dsum\limits_{\Sigma m_{i}=l}\int \dprod\limits_{i=1}^{k}\frac{%
\left( h_{i}/\bar{h}\right) ^{m_{i}}}{m_{i}!}F(\tilde{x})\frac{\partial
^{l+k}\psi }{\partial x_{1}^{m_{1i}+1}...\partial x_{k}^{m_{ki}+1}}(\tilde{x}%
)d\tilde{x}\int K(w)w_{1}^{m_{1}}...w_{k}^{m_{ki}}dw
\end{equation*}
by $(B(h,K),\psi )$ as it represents the value of a linear continuous
functional $B(h,K)$ applied to $\psi .$ The $B(h,K)$ is the leading term in
the generalized bias function for the kernel estimator: 
\begin{eqnarray}
Bias\left( \hat{f}\right) &=&E\hat{f}-f=\overline{h}^{l}B(h,K)+o(\overline{h}%
^{l});  \label{biasexp} \\
\text{where for any }\psi &\in &D_{l+k+1}  \notag \\
\left( E\widehat{f},\psi \right) -\left( f,\psi \right) &=&\overline{h}%
^{l}(B(h,K),\psi )+o(\overline{h}^{l}).  \notag
\end{eqnarray}

The following Theorem gives the limit process for the kernel estimator of
density.

\textbf{Theorem 2. }\textit{For a kernel function }$K$\textit{\ satisfying
Assumption A, if }$\overline{h}\rightarrow 0$\textit{\ and }$\bar{h}%
^{2l}n=O(1)$ \textit{as }$n\rightarrow \infty $\textit{\ the sequence of
generalized random processes }$n^{\frac{1}{2}}\left( \widehat{f}-f-\bar{h}%
^{l}B(h,K)\right) $\textit{\ converges to a generalized Gaussian process
with mean functional zero and covariance functional }$C$\textit{\ which for
any }$\psi _{1},\psi _{2}\in D_{l+k}$\textit{\ \ provides} 
\begin{equation}
(C,(\psi _{1},\psi _{2}))=E\left( \left[ \psi _{1}(x)-E\psi _{1}\left(
x\right) \right] [\psi _{2}(x)-E\psi _{2}x)]\right) =cov\left( \psi
_{1},\psi _{2}\right) .  \label{covar}
\end{equation}%
\textit{If }$n\bar{h}^{2l}\rightarrow 0,$\textit{\ then }$\hat{f}-f$\textit{%
\ converges at the parametric rate }$\sqrt{n}$\textit{\ to a generalized
zero mean Gaussian process with covariance functional }$C$\textit{\ in }$%
\left( \ref{covar}\right) .$

Proof. See appendix.

The condition on the bandwidth that makes it possible to eliminate the bias
asymptotically is less stringent than in the usual topologies and also than
that originally stated in Zinde-Walsh (2008). Under this requirement on the
bandwidth convergence is actually at a parametric rate and the limit
covariance does not involve the kernel function.

\section{\protect\bigskip Distribution function conditional on some
variables and conditional density in the space of generalized functions}

Conditioning is an awkward operation as discussed e.g. in Chang and Pollard
(1997). Here the question posed is limited to conditioning on a variable or
vector in a joint distribution, that is given a joint distribution function $%
F_{x,y}(.,.)$ on $R^{d_{x}}\times R^{d_{y}}$ define a (generalized) function 
$F_{y|x}(.,.)$ that represents the conditional distribution of $y$ given $x.$
A problem associated with such conditioning is that the conditional
distribution function may not exist for every point $x.$

Denote by $F_{x},$ $F_{y}$ the marginal distribution functions of $x,$ $y,$
correspondingly.

Consider limits of ratios to define conditioning:%
\begin{equation}
F_{y|x}=\underset{\Delta \rightarrow 0}{\lim }\frac{F_{x,y}\left( x+\Delta
,y\right) -F_{x,y}\left( x,y\right) }{F_{x}\left( x+\Delta \right)
-F_{x,y}\left( x\right) }.  \label{conddist}
\end{equation}%
As discussed is numerous papers there is a problem defining such a limit
(e.g. Pfazagle, 1979); here it will be demonstrated that the limit exists in
a particular space of generalized functions. Assume that the distribution
function $F_{x}$ is continuous; continuity of this distribution of course
does not preclude singularity.

\textbf{Assumption 2.} \textit{The marginal distribution function }$%
F_{x}\left( x\right) $\textit{\ is continuous on }$R^{d_{x}}.$

Note that although support of the random $y$ belongs to $R^{d_{y}}$ it could
be a discrete set of points, thus we do not restrict $y$ to be continuously
distributed.

Consider the copula function (Sklar, 1973): $C_{F_{x},F_{y}}(a,b)$ on $%
W=\left( 0,1\right) ^{2}$ that is identical to the joint distribution
function, that is for the mapping $M:R^{d_{x}}\times R^{d_{y}}\rightarrow W$
defined by $\left\{ x,y\right\} \rightarrow \left\{
F_{x}(x),F_{y}(y)\right\} $ we get the corresponding mapping $M^{\ast
}(F_{x,y}\left( x,y\right) )=C_{M\left( x,y\right) }\left( M(x,y)\right) $
with 
\begin{equation*}
C_{M\left( x,y\right) }\left( M(x,y)\right)
=C_{F_{x},F_{y}}(F_{x}(x),F_{y}\left( y\right) )=F_{x,y}\left( x,y\right) .
\end{equation*}%
Thus $\left( \ref{conddist}\right) $ is equivalent to%
\begin{equation*}
F_{y|x}=\underset{\Delta \rightarrow 0}{\lim }\frac{C_{F_{x},F_{y}}\left(
F_{x}(x+\Delta ),F_{y}(y)\right) -C_{F_{x},F_{y}}(F_{x}(x),F_{y}\left(
y\right) )}{F_{x}\left( x+\Delta \right) -F_{x,y}\left( x\right) };
\end{equation*}%
denote $F_{x}(x+\Delta )-F_{x}(x)$ by $\tilde{\Delta},$ then by Assumption
2, continuity of $F_{x},$ $\Delta \rightarrow 0$ implies $\tilde{\Delta}%
\rightarrow 0$ thus the limit is equivalent to 
\begin{equation*}
\underset{\tilde{\Delta}\rightarrow 0}{\lim }\frac{C_{F_{x},F_{y}}\left( a+%
\tilde{\Delta},b\right) -C_{F_{x},F_{y}}(a,b)}{\tilde{\Delta}}.
\end{equation*}

Since with respect to its second argument the copula function and the limit
are ordinary functions we concentrate on being able to define the
generalized derivative with respect to the first argument. In particular,
for any $\psi \in D\left( W\right) ,$ given the second argument the value of
the functional $\left( \left( C_{F_{x},F_{y}}\right) _{1}^{\prime },\psi
\right) =-\left( C_{F_{x},F_{y}},\psi ^{\prime }\right) .$ This implies that
we can define the value of the functional $F_{y|x}$ on $D\left( W\right) $
by 
\begin{equation}
\left( F_{y|x},\psi \right) =-\left( C_{F_{x},F_{y}},\psi ^{\prime }\right)
=-\int F_{x,y}(x,y)\psi ^{\prime }\left( F_{x}\left( x\right) \right)
dF_{x}\left( x\right) .  \label{conddistrfl}
\end{equation}%
Thus we can define the conditional distribution $F_{y|x}$ as a generalized
function in the space $D^{\ast }\left( W\right) .$

When $d_{x}=1$ this is an exhaustive representation. When $d_{x}>1$ it may
be advantageous to consider a derivative with respect to a $d_{x}-$%
dimensional argument. Consider the conditioning vector, $x\,,$
component-wise, and consider the multivariate copula function, $%
C_{F_{x_{1}},...,F_{x_{d}},F_{y}}\left( F_{x_{1}},...,F_{x_{d}},F_{y}\right)
;$ to simplify notation we drop the subscript to denote it simply by $C.$
Then by a similar argument for any $\psi \in D\left( W\right) $ where $%
W=\left( 0,1\right) ^{d_{x}}$ we obtain $\left( F_{y|x},\psi \right)
=(-1)^{d_{x}}\left( C,\partial ^{d_{x}}\psi \right) =$ 
\begin{equation}
(-1)^{d_{x}}\int ...\int F_{x,y}(x,y)\partial ^{d_{x}}\psi \left(
F_{x_{1}}\left( x_{1}\right) ,...,F_{x_{d_{x}}}(x_{d_{x}}\right)
)dF_{x_{1}}\left( x_{1}\right) ...dF_{x_{d_{x}}}(x_{d_{x}}).
\label{conddistrmulti}
\end{equation}

\textbf{Remark 1.} Similarly to Corollary 1, the generalized function $%
F_{y|x}$ can be extended as a linear continuous functional from being
defined on the space $D\left( W\right) $ of infinitely differentiable
functions to a linear continuous functional defined by $\left( \ref%
{conddistrfl}\right) $ on any space $D_{k}\left( W\right) $ with $k\geq 1$
and for $\left( \ref{conddistmulti}\right) $ to $D_{k}\left( W\right) $ for
the corresponding $W$ and $k\geq d_{x}.$

\textbf{Remark 2.} If the function $C$ were suitably differentiable the
functional $\left( F_{y|x},\psi \right) $ would be defined for any
continuous $\psi $ with bounded support, that is on the space $D_{0}\left(
W\right) $ by $\left( \partial ^{d_{x}}C\left( ...,.\right) ,\psi \right) :$%
\begin{equation}
\left( F_{y|x},\psi \right) =\int ...\int \partial
^{d_{x}}C(F_{x_{1}},...F_{x_{d_{x}}},F_{y})\psi
(F_{x_{1}},...F_{x_{d_{x}}})dF_{x_{1}}...dF_{x_{d_{x}}}.
\label{conddistrond0}
\end{equation}

In the $y$ argument the conditional distribution is an ordinary function so
here $y$ is considered just as a parameter of the generalized function.
However, the definition of $F_{y|x}$ in $\left( \ref{conddistrfl}\right) $\
can be extended to a functional for functions defined on the product space;
for any $\psi _{x,y}=\psi _{x}(x_{1},...x_{d_{x}})\psi
_{y}(y_{1},...,y_{d_{y}})\in D(\left( 0,1\right) ^{d_{x}})\times
D(R^{d_{y}}) $ define the value of the functional by $\left( F_{y|x},\psi
_{x,y}\right) =$ 
\begin{equation*}
(-1)^{d_{x}}\int ...\int F\left( x,y\right) \partial ^{d_{x}}\psi
_{x}(F_{x_{1}},...F_{x_{d_{x}}})\psi
_{y}(y_{1},...,y_{d_{y}})dF_{x_{1}}...dF_{x_{d_{x}}}dy_{1}...dy_{d_{y}}.
\end{equation*}%
To define conditional density $f_{y|x}$\ as a generalized function one would
have $\left( f_{y|x},\psi _{x,y}\right) =$ 
\begin{equation}
(-1)^{d_{x}+d_{y}}\int ...\int F_{x,y}(x,y)\partial ^{d_{x}}\psi
_{x}(F_{x_{1}},...F_{x_{d_{x}}})\partial ^{d_{y}}\psi
_{y}(y_{1},...,y_{d_{y}})dF_{x_{1}}\left( x_{1}\right)
...dF_{x_{d_{x}}}(x_{d_{x}})dy_{1}...dy_{d_{y}}.  \label{conddensfl}
\end{equation}

In general, the conditional distribution and conditional density depend on
the conditioning variables, $x,$ via the marginals, $F_{x}$; considering
generalized functions makes this explicit.

There are cases when the conditional distribution and conditional density
are defined on the Euclidean space $R^{d_{x}}.$ This is possible if the
distribution function $F_{x}$ is strictly monotone in each argument; then
the corresponding generalized density function is positive, moreover, since
a monotone function is a.e. differentiable, $\partial ^{d_{x}}F_{x,y}(x,y)$
and $f_{x}(x)=\partial ^{d_{x}}F_{x}(x)$ exist a.e. and $f_{x}(x)>0.\,$\
When the density $f_{x}$ is a continuous function the conditional
distribution can be represented as a functional on a function space on $%
R^{d_{x}}$ that can be derived from the general representation above in $%
D^{\ast }\left( W\right) .$

Indeed, any distribution function, $F\left( x,y\right) ,$\ where we focus on
the argument $x,$ via the copula representation can be considered as a
functional on $D\left( W\right) .$ Let $\Phi $ denote the class of such
distribution functions, then $\Phi \subset D^{\ast }\left( W\right) .$
Moreover the representation $\left( \ref{conddistrmulti}\right) $
demonstrated that any conditional distribution $F_{|x}\left( x,y\right) $
also defines a linear continuous functional on $D\left( W\right) .$ Denoting
by $\Phi _{|x}$ the class of conditional distributions we thus have shown
that $\Phi _{|x}\subset D^{\ast }\left( W\right) .$ By the remark, we can
relax the differentiability conditions and consider $\Phi _{|x}\subset
D_{k}^{\ast }\left( W\right) ;$ when the distribution function is
differentiable in $x,$ we set $k=0.$ On the other hand, then a continuous
density function, $f_{x}>0$ exists and the conditional distribution can be
represented by an ordinary function $\frac{\partial ^{d_{x}}F_{x,y}(x,y)}{%
f_{x}(x)};$ denote by $\Phi _{c}$ the class of distributions that are
continuously differentiable in $x$ with $f_{x}>0$ on $R^{d_{x}},$ and by $%
\Phi _{c|x}$ the class of corresponding conditional distributions. Then $%
\Phi _{c}\subset D_{0}^{\ast }(R^{d_{x}})$ and as well $\Phi _{c|x}\subset
D_{0}^{\ast }\left( R^{d_{x}}\right) ,$ where the space $D_{0}\left(
R^{d_{x}}\right) $ is the space of continuous functions with bounded support
in $R^{d_{x}}$. Since $\Phi _{c|x}\subset \Phi _{|x},$ any conditional
distribution that exists in the ordinary sense and thus is in $\Phi _{c|x},$
has two representations: one as a functional on $D_{0}\left( W\right) $
defined above and the second as a functional on $D_{0}\left(
R^{d_{x}}\right) $ that provides for any $\tilde{\psi}\in D_{0}\left(
R^{d_{x}}\right) $%
\begin{equation}
\left( F_{y|x},\tilde{\psi}\right) =\int ...\int \frac{\partial
^{d_{x}}F_{x,y}(x,y)}{f_{x}(x)}\tilde{\psi}(x)dx_{1}...dx_{d_{x}}.
\label{conddistrwithdens}
\end{equation}%
The following lemma shows that the two representations are compatible and
each can be easily obtained from the other.

\textbf{Lemma.}\textit{\ Suppose that }$F_{x,y}\in \Phi _{c}$\textit{. Then
the value of the functional given by }$\left( \ref{conddistrond0}\right) $%
\textit{\ for }$\psi \in D_{0}\left( 0,1\right) ^{d_{x}}$\textit{\ is the
same as the value of the functional given by }$\left( \ref{conddistrwithdens}%
\right) $\textit{\ for }$\tilde{\psi}\left( x\right) =f_{x}(x)\psi \left(
F(x)\right) \in D_{0}\left( R^{d_{x}}\right) ;$\textit{\ and vice versa:
given }$\left( \ref{conddistrwithdens}\right) $\textit{\ the value of }$%
\left( \ref{conddistrond0}\right) $\textit{\ for }$\psi \left(
F_{x_{1}},...,F_{x_{d_{x}}}\right) =\frac{\tilde{\psi}(x_{1},...,x_{d_{x}})}{%
f_{x}(x_{1},...,x_{d_{x}})},$\textit{\ where }$x_{i}$\textit{\ is uniquely
determined by the value of }$F_{x_{i}}$\textit{: }$%
x_{i}=F_{x_{i}}^{-1}(F_{x_{i}}\left( x_{i}\right) ),$\textit{\ is the same.}

Proof. For any $\psi \in D\left( 0,1\right) ^{d_{x}}$ define $\tilde{\psi}$
on $R^{d_{x}}$ by $\tilde{\psi}\left( x\right) =f_{x}(x)\psi \left(
F(x)\right) ,$ then $\left( F_{y|x},\psi \right) $ defined by $\left( \ref%
{conddistrmulti}\right) $ by differentiability of $F_{x,y}\,$\ in $x$ is
equal to 
\begin{equation*}
\left( F_{y|x},\tilde{\psi}\right) =\int ...\int \frac{\partial
^{d_{x}}F_{x,y}(x,y)}{f_{x}(x)}\tilde{\psi}(x)dx_{1}...dx_{d_{x}}.
\end{equation*}%
Denote by $z_{i}$ the value $F_{x_{i}}(x),$ $i=1,...,d_{x};$ then (for
clarity we subscript the operator $\partial $ by the variable(s) with
respect to which we differentiate): 
\begin{equation*}
\partial _{z}^{d_{x}}F_{x,y}\left(
F_{x_{1}}^{-1}(z_{1}),...,F_{x_{d_{x}}}^{-1}\left( z_{d_{x}}\right)
,y\right) f_{x}\left( x\right) =\partial _{x}^{d_{x}}F_{x,y}(x,y).
\end{equation*}%
The r.h.s. of $\left( \ref{conddistrmulti}\right) $ provides%
\begin{eqnarray*}
&&(-1)^{d_{x}}\int ...\int F_{x,y}\left(
F_{x_{1}}^{-1}(z_{1}),...,F_{x_{d_{x}}}^{-1}\left( z_{d_{x}}\right)
,y\right) \partial _{z}^{d_{x}}\psi \left( z_{1},...,z_{d_{x}}\right)
)dz_{1}...dz_{d_{x}} \\
&=&\int ...\int \partial _{z}^{d_{x}}F_{x,y}\left(
F_{x_{1}}^{-1}(z_{1}),...,F_{x_{d_{x}}}^{-1}\left( z_{d_{x}}\right)
,y\right) \psi (z_{1},...z_{d_{x}})dz_{1}...dz_{d_{x}} \\
&=&\int ...\int \frac{\partial _{x}^{d_{x}}F_{x,y}(x,y)}{f_{x}\left(
x\right) }\psi \left( F_{x_{1}}(x_{1}),...,F_{x_{d_{x}}}(x_{d_{x}}\right)
)f_{x}\left( x\right) dx_{1}...dx_{d_{x}}, \\
&&\text{and writing this in more concise notation} \\
&=&\int \frac{\partial ^{d_{x}}F_{x,y}(x,y)}{f_{x}(x)}\psi
(F(x))f_{x}(x)dx=\int \frac{\partial ^{d_{x}}F_{x,y}(x,y)}{f_{x}(x)}\tilde{%
\psi}(x)dx.
\end{eqnarray*}

Since $f_{x}$ is continuous, then $\tilde{\psi}(x)=\psi (F(x))f(x)$ is
continuous on $R^{d_{x}}.$

For an arbitrary $\tilde{\psi}\in D_{0}\left( R^{d_{x}}\right) $ consider 
\begin{equation*}
\left( F_{y|x},\tilde{\psi}\right) =\int \frac{\partial _{x}F_{x,y}(x,y)}{%
f_{x}(x)}\tilde{\psi}(x)dx_{1}...dx_{d_{x}}.
\end{equation*}

Do the transformation, then\ 
\begin{equation*}
\left( F_{y|x},\tilde{\psi}\right) =\int \partial
_{z}F_{x,y}(F_{x}^{-1}(z),y)\frac{\tilde{\psi}(F_{x}^{-1}(z))}{f_{x}\left(
F_{x}^{-1}(z)\right) }dz.
\end{equation*}%
Define a continuous function $\psi \left( F_{x_{1}},...,F_{x_{d_{x}}}\right)
=\frac{\tilde{\psi}(x_{1},...,x_{d_{x}})}{f_{x}(x_{1},...,x_{d_{x}})}$ on $%
\left( 0,1\right) ^{d_{x}},$ then this equals $\left( \ref{conddistrond0}%
\right) .$

$\blacksquare $

Suppose now that $F_{x}$ is absolutely continuous with continuous density
function, $f_{y|x}$; then the support of the density function is an open set
in $R^{d_{x}},$ $S_{y|x}.$ The Lemma applies by considering $\tilde{\psi}%
\left( x\right) =f_{x}(x)\psi \left( F(x)\right) \in D_{0}\left(
S_{y|x}\right) $ in place of $D_{0}\left( R^{d_{x}}\right) .$

\section{Limit properties of kernel estimators of conditional distribution
in generalized functions}

Consider the usual kernel estimator of conditional distribution; typically
its limit properties are available under smoothness conditions on the
distribution (see, e.g. Li and Racine, 2007). Here the estimator is examined
in the space of generalized functions without any restrictions placed on the
distribution beyond Assumption 2 (continuity of $F_{x}).$

Recall the usual kernel estimator of conditional distribution:%
\begin{eqnarray}
\hat{F}_{y|x}\left( x,y\right) &=&\frac{\Sigma \bar{G}\left( \frac{y-y_{i}}{%
h_{y}}\right) K\left( \frac{x_{i}-x}{h}\right) }{\Sigma K\left( \frac{x_{i}-x%
}{h}\right) }  \label{conddistrest} \\
&=&\frac{\frac{1}{n}\Sigma \bar{G}\left( \frac{y-y_{i}}{h_{y}}\right) \frac{1%
}{h^{d_{x}}}K\left( \frac{x_{i}-x}{h}\right) }{\hat{f}_{x}(x)},
\end{eqnarray}%
where $\bar{G}$ is the integral of a kernel function $G$ similar to $K$ that
satisfies Assumption 1 on $R^{d_{y}}$ and $K$ satisfies Assumption 1 on $%
R^{d_{x}}.$ Sometimes $\bar{G}$ is assumed to be the indicator function $%
I(w>0).$

To simplify exposition we assume that each component of vector $x$ is
associated with the same (scalar) bandwidth parameter $h;$ it is not
difficult to generalize to the case of distinct bandwidths.

\textbf{Theorem 3.} \textit{Suppose that Assumption 1 on the kernel }$K$%
\textit{\ and either a similar assumption for }$G$\textit{\ holds, or }$\bar{%
G}$ \textit{is the indicator function, the bandwidth parameter }$%
h=cn^{-\alpha },$\textit{\ where }$\alpha <\frac{1}{4}$\textit{\ and
Assumption 2 holds. Then for a random sample }$\left\{ \left(
x_{i},y_{i}\right) \right\} _{i=1}^{n}$ \textit{the estimator }$\hat{F}%
_{y|x}\left( x,y\right) $\textit{\ as a generalized random function on }$%
D\left( W\right) $\textit{\ converges to the conditional distribution
generalized function }$F_{y|x}$\textit{\ defined by }$\left( \ref%
{conddistrfl}\right) $\textit{\ at the rate }$n^{-\frac{1}{2}};$ \textit{the
limit process for }$\sqrt{n}(\hat{F}_{y|x}-F_{y|x})$\textit{\ on }$D\left(
W\right) $\textit{\ is given by a }$\psi \in D\left( W\right) $ \textit{%
indexed random functional, }$Q_{y|x}$ with $(Q_{y|x},\psi )=$\textit{\ \ }%
\begin{equation*}
\left( -1\right) ^{d_{x}}[\int F_{xy}(\partial ^{d_{x}}\partial ^{d_{x}}\psi
)(F_{x})U_{x}dF_{x}+\int F_{xy}(\partial ^{d_{x}}\psi )(F_{x})dU_{x}+\int
(\partial ^{d_{x}}\psi )(F_{x})U_{xy}dF_{x}],
\end{equation*}%
\textit{where }$U_{x},U_{xy}$\textit{\ are Brownian bridge processes with
dimension }$d_{x},d_{y}+d_{x},$\textit{\ correspondingly; as a generalized
random process the limit process }$Q_{y|x}$\textit{\ of }$\sqrt{n}(\hat{F}%
_{y|x}-F_{y|x})$\textit{\ is Gaussian with mean functional zero and
covariance bilinear functional }$C,$\textit{\ given for any }$\psi _{1},\psi
_{2}$ \textit{by }%
\begin{equation*}
\left( C,(\psi _{1},\psi _{2}\right) )=cov[(Q_{y|x},\psi _{1}),(Q_{y|x},\psi
_{2}).
\end{equation*}

Proof. See Appendix.

This result is general in that the root-n convergence holds here regardless
of whether the marginal density exists. If it does exist the result could be
restated for conditional distribution as a generalized function on $%
D_{0}\left( R^{d_{x}}\right) $ by $\left( \ref{conddistrwithdens}\right) .$

\textbf{Remark 3.} Sometimes for a singular distribution the kernel
estimator $\hat{f}_{x}\left( x\right) $ diverges at a specific rate, as e.g.
in Lu $\left( 1999\right) $ where at points $x$ in support of density $\hat{f%
}_{x}\left( x\right) =h^{d-1}b+o_{p}\left( h^{d-1}\right) $ with some $b>0$
and $d=\frac{\ln 2}{\ln 3}<1$. In the univariate case this is discussed in
Example 5 in Zinde-Walsh (2008), where for the Cantor distribution it is
noted that though $\hat{f}_{x}(x)$ may diverge, $h^{1-d}\hat{f}_{x}\left(
x\right) $ is bounded and bounded away from zero. Then, even though the
limit density does not exist by rescaling it is possible to establish the
convergence rate of the estimator of the conditional distribution as a
functional on $D_{0}\left( R^{d_{x}}\right) ;$ the rate is $n^{-\frac{1}{2}%
}h^{1-d}$ and is faster than the root-n rate.

\section{Conditional moments}

Consider now a conditional moment of a function $g\left( y\right) ,$ of $%
y\in R^{d_{y}}:$ $E_{y|x}g(y)=m\left( x\right) ,$ with $m\left( x\right) $
measurable with respect to $F_{x}.$

When the conditional density function exists in $L_{1}$ we write $m\left(
x\right) =\int g(y)f_{y|x}(x,y)dy$ (assuming that the integral exists). As a
generalized function (in $x)$ $m\left( x\right) $ can be presented on the
space $D\left( W\right) ;$ $W=\left( 0,1\right) ^{d_{x}}$ by the value of
the functional for $\psi :$%
\begin{equation*}
\left( m,\psi \right) =\int m\left( x\right) \psi (F(x))dF(x)=\int \left[
\int g(y)f_{y|x}(x,y)dy\right] \psi (F(x))dF(x).
\end{equation*}

To give meaning to $\left( m,\psi \right) $ regardless of the existence of
the conditional density as a function, $\int g(y)f_{y|x}(x,y)dy$ needs to be
characterized as a generalized function on $D\left( W\right) .$ To make this
possible for an arbitrary distribution on $\left( x,y\right) $ that
satisfies Assumption 2 the class of functions $g$ is restricted.

\textbf{Assumption 3. }\textit{The function }$g$ \textit{is continuously
differentiable with respect to the differentiation operator }$\partial
^{d_{y}}$\textit{.}

Any polynomial function satisfies Assumption 3, and thus conditional mean of 
$y,$ or conditional variance (if they exist) can be considered. If the
function were not to satisfy the differentiability assumption, the class of
distributions would need to be correspondingly restricted.

Consider $D\left( R^{d_{y}}\right) $ and a locally finite partition of unity
on $R^{d_{y}}$ by a set of suitable functions, "bump" functions from $%
D\left( R^{d_{y}}\right) :\left\{ \psi _{\nu }\right\} ,$ where $\psi _{\nu
}\in D\left( R^{d_{y}}\right) ,$ $\psi \geq 0$ and $\Sigma _{\nu }\psi _{\nu
}\left( y\right) \equiv 1;$ also any $y$ can belong to support of only a
finite number of $\psi _{v}.$ See e.g. Gel'fand and Shilov, 1964, v.1, p.142
for a construction.

Then define $\left( gf_{y|x},\psi _{\nu }\right) =\int g(y)f_{y|x}(x,y)\psi
_{v}(y)dy;$ under Assumption 3 this expression is (as usual integrating by
parts and using boundedness of support of $\psi _{v}$):%
\begin{equation}
\int g(y)f_{y|x}(x,y)\psi _{v}(y)dy=\left( -1\right) ^{d_{y}}\int
F_{y|x}\left( x,y\right) \partial ^{d_{y}}\left( g\left( y\right) \psi
_{v}\left( y\right) \right) dy.  \label{eq}
\end{equation}%
This expression represents a generalized function on $D\left( W\right) $
given for any $\psi \in D\left( W\right) $ by%
\begin{eqnarray*}
&&(\int g(y)f_{y|x}(x,y)\psi _{v}(y)dy,\psi ) \\
&=&\left( -1\right) ^{d_{y}}\int \int F_{y|x}\left( x,y\right) \partial
^{d_{y}}\left( g\left( y\right) \psi _{v}\left( y\right) \right) dy\psi
(F(x))dF(x) \\
&=&\left( -1\right) ^{d_{y}+d_{x}}\int \int F_{x,y}\left( x,y\right)
\partial ^{d_{y}}\left( g\left( y\right) \psi _{v}\left( y\right) \right)
dy(\partial ^{d_{x}}\psi )(F(x))dF(x).
\end{eqnarray*}%
Because the supports of $\psi _{v}$ and of $\psi $ are bounded and the
function being integrated is bounded, the integral exists.

\textbf{Assumption 4.} (Existence of conditional moment). For a partition of
unity, $\left\{ \psi _{\nu }\right\} ,$ the sum 
\begin{equation}
\Sigma _{v}(\int g(y)f_{y|x}(x,y)\psi _{v}(y)dy,\psi )  \label{sum}
\end{equation}%
converges.

Then $\left( \ref{sum}\right) $ represents $\left( m\left( x\right) ,\psi
\right) $ for the generalized function, $m\left( x\right) =\Sigma _{v}\int
g(y)f_{y|x}(x,y)\psi _{v}(y)dy,$ on $D\left( W\right) .$

Thus 
\begin{equation*}
m\left( x\right) =\int g(y)f_{y|x}(x,y)dy=\Sigma _{v}\left( gf_{y|x},\psi
_{\nu }\right) ,
\end{equation*}%
where the sum converges.

Then 
\begin{eqnarray*}
\Sigma _{v}\int g(y)f_{y|x}(x,y)\psi _{v}(y)dy &=&\int
g(y)f_{y|x}(x,y)\Sigma _{v}\psi _{v}\left( y\right) dy \\
&=&\int g(y)f_{y|x}(x,y)dy,
\end{eqnarray*}%
in other words interchanging the order of integration and summation is
permitted for the terms on the left-hand side of $\left( \ref{eq}\right) $
under Assumption 4. However, this is not the case for terms on the
right-hand side of $\left( \ref{eq}\right) .$ For example, if $g\left(
y\right) =y,$ we have $\partial ^{d_{y}}\left( g\left( y\right) \psi
_{v}\left( y\right) \right) =y\psi _{v}^{\prime }+\psi _{v},$ and $\Sigma
_{v}\left( \partial ^{d_{y}}\left( g\left( y\right) \psi _{v}\left( y\right)
\right) \right) =1,$ but $\int F_{y|x}\left( x,y\right) dy$ may not exist.

Thus $\left( gf_{y|x},\psi \psi _{\nu }\right) =$ 
\begin{equation}
(-1)^{d_{x}+d_{y}}\int ...\int F_{x,y}(x,y)\partial ^{d_{x}}\psi \left(
F_{x}\left( x\right) )\partial ^{d_{y}}[g(y)\psi
_{v}(y_{1},...,y_{d_{y}}\right) ]dF_{x}(x)dy_{1}...dy_{d_{y}}.
\label{condmomv}
\end{equation}

Then the conditional moment $m$ as a generalized function on $D\left(
W\right) $ is given by $\left( m,\psi \right) =$%
\begin{equation}
\Sigma _{v}(-1)^{d_{x}+d_{y}}\int ...\int F_{x,y}(x,y)\partial ^{d_{x}}\psi
\left( F_{x}\left( x\right) )\partial ^{d_{y}}[g(y)\psi
_{v}(y_{1},...,y_{d_{y}}\right) ]dF_{x}(x)dy_{1}...dy_{d_{y}}
\label{condmom}
\end{equation}%
with any $\left\{ \psi _{v}\right\} $ representing a partition of unity on $%
R^{d_{y}}$ by functions from $D\left( R^{d_{y}}\right) .$

\section{Limit properties of kernel estimators of conditional mean function.}

Suppose that with $d_{y}=1$ the conditional mean function $m\left( x\right)
=E_{y|x}y$ exists; by $\left( \ref{condmom}\right) $ it then can be
represented as%
\begin{eqnarray}
&&\left( m,\psi \right)  \label{mkci} \\
&=&\Sigma _{v}(-1)^{d_{x}+1}\int ...\int F_{x,y}(x,y)\partial ^{d_{x}}\psi
\left( F_{x}\left( x\right) )[y\psi _{v}^{\prime }(y\right) +\psi _{v}\left(
y\right) ]dF_{x}\left( x\right) dy_{1}...dy_{d_{y}}.  \notag
\end{eqnarray}

Consider the usual kernel estimator 
\begin{equation*}
\hat{m}(x)=\frac{\Sigma y_{i}K\left( \frac{x_{i}-x}{h}\right) }{\Sigma
K\left( \frac{x_{j}-x}{h}\right) },
\end{equation*}%
that can also be represented as%
\begin{equation*}
\frac{\int y\hat{f}_{x,y}(x,y)dy}{\hat{f}_{x}(x)}=\frac{\Sigma _{v}\int y%
\hat{f}_{x,y}(x,y)\psi _{v}(y)dy}{\hat{f}_{x}(x)}.
\end{equation*}%
Then for any continuously differentiable $\tilde{\psi}(x)$ 
\begin{eqnarray*}
\left( \hat{m},\tilde{\psi}\right) &=&\tint \frac{\Sigma _{v}\int y\hat{f}%
_{x,y}(x,y)\psi _{v}(y)dy}{\hat{f}_{x}(x)}\tilde{\psi}(x)dx \\
&=&-\Sigma _{v}\int \frac{\int \partial ^{d_{x}}\hat{F}_{x,y}(x,y)\frac{d}{dy%
}[y\psi _{v}(y)]dy}{\hat{f}_{x}(x)}\tilde{\psi}(x)dx \\
&=&-\Sigma _{v}(\hat{m},\tilde{\psi}\psi _{v}).
\end{eqnarray*}
Consider $\psi $ and $\tilde{\psi}=\psi \hat{f}$ ; by the Lemma $(\hat{m},%
\tilde{\psi}\psi _{v})=$ 
\begin{eqnarray}
&&\left( -1\right) ^{d_{x}+1}\int \int \hat{F}_{x,y}\left( x,y\right)
\partial ^{d_{x}}\psi \left( \hat{F}_{x}(x)\right) \frac{d}{dy}[y\psi
_{v}(y)]d\left( \hat{F}_{x}(x)\right) dy  \label{meanfunc} \\
&=&\left( -1\right) ^{d_{x}+1}\int \int \hat{F}_{x,y}\left( x,y\right)
\partial ^{d_{x}}\psi \left( \hat{F}_{x}(x)\right) [y\psi _{v}^{\prime
}(y)+\psi _{v}\left( y\right) ]d\left( \hat{F}_{x}(x)\right) dy.  \notag
\end{eqnarray}

\bigskip \textbf{Assumption 5. }The conditional variance $\sigma ^{2}\left(
x\right) =E_{y|x}y^{2}$ defines a generalized function on $D\left( W\right) $%
.

Assumption 5 implies that for any $\psi \in D\left( W\right) $ the value of
the functional $\left( \sigma ^{2},\psi \right) =\int \sigma ^{2}\left(
x\right) \psi \left( F_{x}(x)\right) dF_{x}(x)$ is always bounded; this is
reqired to bound the variance for the limit process. By $\left( \ref{condmom}%
\right) $ for a partition of unity, $\left\{ \psi _{v}\right\} $%
\begin{equation*}
\left( \sigma ^{2},\psi \right) =\Sigma _{v}\left( -1\right) ^{d_{x}+1}\int
\int F_{x,y}(x,y)\partial ^{d_{x}}\psi (F_{x}\left( x\right) )\left(
y^{2}\psi _{v}\left( y\right) \right) ^{\prime }dF_{x}\left( x\right) dy.
\end{equation*}

\textbf{Theorem 4. }\textit{Suppose that Assumptions 1-5 hold, the bandwidth
parameter }$h=cn^{-\alpha },$\textit{\ where }$\alpha <\frac{1}{4}$\textit{.
Then the estimator }$\hat{m}(x)$\textit{\ for a random sample }$\left\{
\left( x_{i},y_{i}\right) \right\} _{i=1}^{n}$ \textit{as a generalized
random function on }$D\left( W\right) $\textit{\ converges at the rate }$n^{-%
\frac{1}{2}}$ \textit{to the generalized function }$m$\textit{\ that
provides }$\left( \ref{mkci}\right) ;$ \textit{the limit process for }$\sqrt{%
n}(\hat{m}-m)$\textit{\ on }$D\left( W\right) $\textit{\ is given by a }$%
\psi \in D\left( W\right) $\textit{\ indexed random \ functional }$Q_{m}$
with $(Q_{m},\psi )=$%
\begin{eqnarray*}
&&\Sigma _{v}(-1)^{d_{x}+1}\int ...\{\int U_{x,y}\partial ^{d_{x}}\psi
(F_{x}\left( x\right) )dF_{x}\left( x\right) \\
&&+\int F_{x,y}\left( x,y\right) \left( \partial ^{d_{x}}\right) ^{2}\psi
(F_{x}\left( x\right) )U_{x}dF_{x}\left( x\right) \\
&&+\int F_{x,y}\left( x,y\right) \partial ^{d_{x}}\psi (F_{x}\left( x\right)
)dU_{x}\}[y\psi _{v}^{\prime }(y)+\psi _{v}\left( y\right)
]dy_{1}...dy_{d_{y}},
\end{eqnarray*}%
\textit{where }$U_{x},U_{x,y}$\textit{\ are Brownian bridge processes with
dimension }$d_{x},d_{x}+1,$\textit{\ correspondingly; as a generalized
random process the limit process }$Q_{m}$\textit{\ of }$\sqrt{n}(\hat{m}-m)$%
\textit{\ is Gaussian with mean functional zero and covariance bilinear
functional }$C,$\textit{\ given for any }$\psi _{1},\psi _{2}$ \textit{by }%
\begin{equation*}
\left( C,(\psi _{1},\psi _{2}\right) )=cov[(Q_{m},\psi _{1}),(Q_{m},\psi
_{2}).
\end{equation*}

Proof. See Appendix.

Similarly to the kernel estimator for the conditional distribution the
conditional mean estimator converges at parametric rate as a functional on $%
D\left( W\right) $ for any distribution. When a positive conditioning
density exists it is possible to represent the conditional mean as a
functional on $D\left( R^{d_{x}}\right) ,$ by the same arguments as in the
Lemma. In the case of Remark 3 a similar rescaling provides a faster
convergence rate for the estimator considered as a functional on $D\left(
R^{d_{x}}\right) .$

\section{Conclusion and further questions}

The approach employed here makes it possible to avoid any restrictions when
defining density, conditional distribution and conditional density as well
as conditional moments for a smooth function (e.g. conditional expectation
or second moment).

The usual kernel estimators converge to the limit generalized functions at a
parametric rate; the limit process is provided by a Gaussian process in the
space of generalized functions, that is a Gaussian process indexed by
well-behaved functions from the appropriate spaces.

The results here were based on a random sample of observations to simplify
exposition; extension to stationary ergodic or mixing processes can be
obtained. Further extensions to relax homogeneity and independence are a
subject of future research.

The limit results imply that with a judicial selection of indexing functions
one could use the kernel estimators for inference in very general
situations; this investigation is mostly left for future research.

\section{Appendix.}

\textbf{Proof of Theorem 2.}

Define a generalized function $e_{nhj}$ such that the value of the
functional for $\psi \in G$ is%
\begin{equation*}
(e_{nhj},\psi )=\int \frac{1}{\Pi h_{i}}K(\frac{x-x_{j}}{h})\psi
(x)dx-(f,\psi )
\end{equation*}%
and consider $e_{nh}=\frac{1}{n}\sum_{j=1}^{n}e_{nhj}$; this generalized
function provides $\widehat{f}-f.$

The expectation functional $Ee_{hn}$ gives the generalized bias of the
estimator $\hat{f},$ $Bias\left( \hat{f}\right) ,$ see $\left( \ref{biasexp}%
\right) .$

Next to derive the variance functional consider $T_{lj}=E(e_{nhl},\psi
_{1})(e_{hnj},\psi _{2})).$

For $l\neq j$ by independence 
\begin{eqnarray*}
T_{lj} &=&E(e_{nhl},\psi _{1})(e_{nhj},\psi _{2})=E(e_{nhl},\psi
_{1})E(e_{nhj},\psi _{2}) \\
&=&\left( Bias\left( \hat{f}\right) ,\psi _{1}\right) \left( Bias\left( \hat{%
f}\right) ,\psi _{2}\right) .
\end{eqnarray*}%
For $l=j$ 
\begin{eqnarray*}
T_{jj} &=&E(e_{nhj}(x),\psi _{1})(e_{nhj}(x),\psi _{2}) \\
&=&\int \left[ \int \frac{1}{\Pi h_{i}}K(\frac{x_{j}-x}{h})\psi
_{1}(x)dx-(f,\psi _{1})\right] \times \\
&&\left[ \int \frac{1}{\Pi h_{i}}K(\frac{x_{j}-x}{h})\psi _{2}(x)dx-(f,\psi
_{2})\right] dF(x_{j}) \\
&=&T_{jj}^{1}+T_{jj}^{2},
\end{eqnarray*}%
where 
\begin{equation*}
T_{jj}^{1}=\int \left( \int \frac{1}{\Pi h_{i}}K(\frac{x_{j}-x}{h})\psi
_{1}(x)dx\right) \left( \int \frac{1}{\Pi h_{i}}K(\frac{x_{j}-x}{h})\psi
_{2}(x)dx\right) dF(x_{j})
\end{equation*}%
and $T_{jj}^{2}=$ 
\begin{eqnarray*}
&&-\int [\int \frac{1}{\Pi h_{i}}K(\frac{x_{j}-x}{h})\psi
_{1}(x)dx]dF(x_{j})\times (f,\psi _{1}) \\
&&-\int [\int \frac{1}{\Pi h_{i}}K(\frac{x_{j}-x}{h})\psi
_{2}(x)dx]dF(x_{j})\times (f,\psi _{2}) \\
&&+(f,\psi _{1})\times (f,\psi _{2}).
\end{eqnarray*}

For every vector $h$ and $s=1,2$ 
\begin{equation*}
\int \frac{1}{\Pi h_{i}}K\left( \frac{x_{j}-x}{h}\right) \psi _{s}(x)dx=\int
K\left( w\right) \psi _{s}\left( x_{j}-hw\right) dw.
\end{equation*}

It follows by substituting into $T_{jj}^{2}$ and expanding $\psi _{s}$ that $%
T_{jj}^{2}=-E\psi _{1}\left( x\right) E\psi _{2}\left( x\right) +\bar{h}%
R_{2}.$

Similarly, 
\begin{eqnarray*}
T_{jj}^{1} &=&\int \left( \int \frac{1}{\Pi h_{i}}K\left( \frac{x_{j}-x}{h}%
\right) \psi _{1}(x)dx\right) \left( \int \frac{1}{\Pi h_{i}}K\left( \frac{%
x_{j}-x}{h}\right) \psi _{2}(x)dx\right) dF(x_{j}) \\
&=&\int \left( \int K\left( w\right) \psi _{1}\left( x_{j}-hw\right) dw\int
K\left( w\right) \psi _{2}\left( x_{j}-hw\right) dw\right) dF(x_{j})= \\
&=&\int \left( \int K(w)dw\psi _{1}(x_{j})-\bar{h}\int K\left( w\right) %
\left[ \Sigma _{i=1}^{k}\frac{\partial \psi _{1}}{\partial x_{i}}\left(
x_{j}-h\tilde{w}\right) w_{i}\frac{h_{i}}{\bar{h}}\right] dw\right) \times \\
&&\left( \int K(w)dw\psi _{2}(x_{j})-\bar{h}\int K\left( w\right) \left[
\Sigma _{i=1}^{k}\frac{\partial \psi _{2}}{\partial x_{i}}\left( x_{j}-h%
\tilde{w}\right) w_{i}\frac{h_{i}}{\bar{h}}\right] dw\right) dF(x_{j}) \\
&=&E\psi _{1}\left( x\right) \psi _{2}\left( x\right) +\bar{h}R_{1};
\end{eqnarray*}%
where after the change of variable $\psi _{s}\left( x_{j}-hw\right) $ is
expanded around the point $x_{j}.$ Next we establish that ${\text{ }}%
\left\vert R_{1}\right\vert <\infty ,\left\vert R_{2}\right\vert <\infty .$

Indeed, 
\begin{equation}
\psi _{s}\left( x-hw\right) =\psi _{s}\left( x\right) -\bar{h}\Sigma
_{i=1}^{k}\frac{\partial \psi _{s}}{\partial x_{i}}\left( x-h\tilde{w}%
\right) w_{i}\frac{h_{i}}{\bar{h}},s=1,2,
\end{equation}%
where $\tilde{w}=\alpha w$ for some $0\leq \alpha \leq 1$ and since $%
h_{i}\leq \bar{h}$ and $\left\vert w\right\vert <1$ on support of $K$ 
\begin{equation*}
\left\vert \Sigma _{i=1}^{k}\frac{\partial \psi _{s}}{\partial x_{i}}\left(
x-h\tilde{w}\right) w_{j}\frac{h_{j}}{\bar{h}}\right\vert \leq \left\vert
\Sigma _{i=1}^{k}\frac{\partial \psi _{s}}{\partial x_{i}}\left( x-h\tilde{w}%
\right) \right\vert
\end{equation*}%
holds and the right-hand side is uniformly bounded by some $B_{\psi
_{s}}<\infty $ since $\psi _{s}\in D_{l+k}\left( U\right) .$ Thus 
\begin{equation*}
|R_{1}|\leq B_{\psi _{1}}\sup \psi _{2}+B_{\psi _{2}}\sup \psi _{1}+\bar{h}%
B_{\psi _{1}}B_{\psi _{2}}.
\end{equation*}%
Similarly, $\left\vert R_{2}\right\vert <\infty .$

Combining we get that $T_{jj}=cov\left( \psi _{1},\psi _{2}\right) +O(\bar{h}%
)$ as $\bar{h}\rightarrow 0.$

Consider now 
\begin{eqnarray}
\eta _{nhj} &=&n^{\frac{1}{2}}[e_{nhj}-E(e_{nhj})];  \notag \\
\eta _{nh} &=&\frac{1}{n}\sum \eta _{nhj}.  \label{eta}
\end{eqnarray}%
Note that here $\eta _{nhj}=n^{\frac{1}{2}}(e_{nhj}-Bias\left( \hat{f}%
\right) ).$ This generalized random function has expectation zero. In the
covariance the terms where $l\neq j$ are zero and 
\begin{eqnarray*}
&&n^{-1}E(\eta _{nhj},\psi _{1})(\eta _{nhj},\psi _{2}) \\
&=&T_{jj}+O(\bar{h}),
\end{eqnarray*}
and thus converges to $cov(\psi _{1},\psi _{2}).$

Next (similarly to Zinde-Walsh, 2008) we show that for any set of linearly
independent functions $\psi _{1},...,\psi _{m}\in D$ with $E(\psi
_{l}^{2})>0 $ the joint distribution of the vector 
\begin{equation*}
\overrightarrow{\eta }_{nh}=((\eta _{nh},\psi _{1})...,(\eta _{nh},\psi
_{m}))^{\prime }
\end{equation*}%
converges to a multivariate Gaussian. Define similarly the vector $%
\overrightarrow{\eta }_{nhj}$ with components ($\eta _{nhj},\psi _{l}).$
Denote by $S$ the $m\times m$ matrix with $ts$ component $\left\{ S\right\}
_{ts}=(C,(\psi _{t},\psi _{s}))$ where the functional $C$ is given by (\ref%
{covar})$.$ Denote by $\hat{S}_{n}$ the covariance matrix of $%
\overrightarrow{\eta }_{nhj}.$ By the convergence results for $T_{lj},$ $%
\hat{S}_{n}\rightarrow \Sigma .$ Since the functions $\psi _{1},...,\psi
_{m} $ are linearly independent and $E(\psi _{l}^{2})>0$ the matrix $S$ and
thus $\allowbreak \hat{S}_{n}$ for large enough $n$ is invertible. Define $%
\xi _{nhj}$ to equal $\hat{S}_{n}^{-1/2}\overrightarrow{\eta }_{nhj},$ then $%
\hat{S}_{n}^{-1/2}\overrightarrow{\eta }_{nhj}-S^{-1/2}\overrightarrow{\eta }%
_{nhj}\rightarrow _{p}0.$

Next, consider an $m\times 1$ vector $\lambda $ with $\lambda ^{\prime
}\lambda =1.$ The random variables $\lambda ^{\prime }\xi _{nhj}$ are
independent with expectation 0, $var\sum \lambda ^{\prime }\xi _{nhj}=1;$
they satisfy the Liapunov condition: $\sum E\left\vert \lambda ^{\prime }\xi
_{nhj}\right\vert ^{2+\delta }\rightarrow 0$ for $\delta >0$ since the
kernel function is bounded with finite support. Thus 
\begin{equation*}
\dsum \lambda ^{\prime }\xi _{nhj}\rightarrow _{d}N(0,1)
\end{equation*}%
and by the Cramer-Wold theorem convergence to a limit Gaussian process for $%
\hat{S}_{n}^{-1/2}\overrightarrow{\eta }_{nh}$ and thus for $S^{-1/2}%
\overrightarrow{\eta }_{hn}$ follows. $\blacksquare $

\textbf{Proof of Theorem 3.}

Since for a smooth kernel $\hat{F}(x,y)\in \Phi _{c}$ by the Lemma the value
of the functional for $\psi \in D\left( 0,1\right) ^{d_{x}},$ $(\hat{F}%
_{y|x},\psi )$ is the same as $(\hat{F}_{y|x},\tilde{\psi}),$ with the
latter defined by $\left( \ref{conddistrond0}\right) $ where $\tilde{\psi}=%
\hat{f}_{x}\psi \left( \hat{F}_{x}\right) .$ Thus for any $\psi \in D\left(
0,1\right) :$ 
\begin{eqnarray}
&&(\hat{F}_{y|x},\psi )  \label{conddistrestfl} \\
&=&\left( -1\right) ^{d_{x}}\int \frac{1}{n}\Sigma \bar{G}\left( \frac{%
y-y_{i}}{h_{y}}\right) \bar{K}\left( \frac{x_{i}-x}{h}\right) \partial
^{d_{x}}\psi \left( \Sigma \bar{K}\left( \frac{x_{i}-x}{h}\right) \right)
d\left( \Sigma \bar{K}\left( \frac{x_{i}-x}{h}\right) \right) .
\end{eqnarray}

More concisely it is $(\hat{F}_{y|x},\psi )=$%
\begin{eqnarray*}
&&\left( -1\right) ^{d_{x}}\int \hat{F}_{x,y}\left( x,y\right) \partial
^{d_{x}}\psi \left( \hat{F}_{x}(x)\right) d\left( \hat{F}_{x}(x)\right) \\
&&+\left( -1\right) ^{d_{x}}[\int \tilde{F}_{x,y}\left( x,y\right) \partial
^{d_{x}}\psi \left( \tilde{F}_{x}(x)\right) d\left( \tilde{F}_{x}(x)\right)
-\int \hat{F}_{x,y}\left( x,y\right) \partial ^{d_{x}}\psi \left( \hat{F}%
_{x}(x)\right) d\hat{F}_{x}(x)].
\end{eqnarray*}%
Here "hat" indicates empirical distribution function and "tilde" the kernel
estimated distribution function. By standard arguments the smooth kernel
introduces a bias; by the usual expansions using differentiability of $\psi $
we get that for the second order kernel 
\begin{eqnarray*}
&&\left( -1\right) ^{d_{x}}[\int \tilde{F}_{x,y}\left( x,y\right) \partial
^{d_{x}}\psi \left( \tilde{F}_{x}(x)\right) d\left( \tilde{F}_{x}(x)\right)
-\int \hat{F}_{x,y}\left( x,y\right) \partial ^{d_{x}}\psi \left( \hat{F}%
_{x}(x)\right) d\hat{F}_{x}(x)] \\
&=&O_{p}\left( h^{2}\right) .
\end{eqnarray*}%
Represent $\left( -1\right) ^{d_{x}}\int \hat{F}_{x,y}\left( x,y\right)
\partial ^{d_{x}}\psi \left( \hat{F}_{x}(x)\right) d\left( \hat{F}%
_{x}(x)\right) $ as%
\begin{eqnarray*}
&&\left( -1\right) ^{d_{x}}\{\int F_{x,y}\partial ^{d_{x}}\psi \left(
F_{x}\right) d\left( F_{x}\right) +\int F_{x,y}[(\partial ^{d_{x}}\partial
^{d_{x}}\psi )\left( F_{x}\right) \left( \hat{F}_{x}-F_{x}\right) +r\left( 
\hat{F}_{x}-F_{x}\right) ^{2}]d\left( F_{x}\right) \\
&&+\int F_{x,y}\partial ^{d_{x}}\psi \left( F_{x}\right) d\left( \hat{F}%
_{x}-F_{x}\right) +\int F_{x,y}(\partial ^{d_{x}}\partial ^{d_{x}}\psi
)\left( \tilde{F}_{x}\right) \left( \hat{F}_{x}-F_{x}\right) d\left( \hat{F}%
_{x}-F_{x}\right) \\
&&+\int (\hat{F}_{x,y}-F_{x,y})\partial ^{d_{x}}\psi \left( F_{x}\right)
dF_{x}+\int \left( \hat{F}_{x,y}-F_{x,y}\right) \partial ^{d_{x}}\psi \left(
F_{x}\right) d\left( \hat{F}_{x}-F_{x}\right) \\
&&+\int (\hat{F}_{x,y}-F_{x,y})(\partial ^{d_{x}}\partial ^{d_{x}}\psi
)\left( \tilde{F}_{x}\right) \left( \hat{F}_{x}-F_{x}\right) dF_{x} \\
&&+\int (\hat{F}_{x,y}-F_{x,y})(\partial ^{d_{x}}\partial ^{d_{x}}\psi
)\left( \tilde{F}_{x}\right) \left( \hat{F}_{x}-F_{x}\right) d\left( \hat{F}%
_{x}-F_{x}\right)
\end{eqnarray*}%
where $\tilde{F}_{x}$ represents an intermediate value and takes values in $%
\left( 0,1\right) ^{d_{x}};$ by properties of $\psi \in D\left( W\right) $
the function $(\partial ^{d_{x}}\partial ^{d_{x}}\psi )\left( \tilde{F}%
_{x}\right) $ is bounded. Then $\sqrt{n}\left( \hat{F}_{y|x}-F_{y|x},\psi
\right) $ can be expressed as 
\begin{equation*}
Q_{\psi }\left( \sqrt{n}(\hat{F}_{x}-F_{x}),\sqrt{n}\left( \hat{F}%
_{xy}-F_{xy}\right) \right) +n^{-\frac{1}{2}}R\left( \sqrt{n}(\hat{F}%
_{x}-F_{x}),\sqrt{n}\left( \hat{F}_{xy}-F_{xy}\right) \right) ,
\end{equation*}%
where 
\begin{eqnarray*}
&&Q_{\psi }\left( \sqrt{n}(\hat{F}_{x}-F_{x}),\sqrt{n}\left( \hat{F}%
_{xy}-F_{xy}\right) \right) \\
&=&\int F_{x,y}[(\partial ^{d_{x}}\psi )\left( F_{x}\right) ]d\sqrt{n}\left( 
\hat{F}_{x}-F_{x}\right) +\int \sqrt{n}(\hat{F}_{x,y}-F_{x,y})[(\partial
^{d_{x}}\psi )\left( F_{x}\right) ]dF_{x} \\
&&+\int F_{x,y}[(\partial ^{d_{x}}\partial ^{d_{x}}\psi )\left( F_{x}\right)
]\sqrt{n}\left( \hat{F}_{x}-F_{x}\right) d\left( F_{x}\right)
\end{eqnarray*}%
and $R\left( .,.\right) $ is a bounded function.

Since the limit process of $\sqrt{n}\left( \hat{F}_{.}-F_{.}\right) $ is $%
U_{.},$ a Brownian bridge, and the function $Q_{\psi }$ is continuous in its
arguments, by Donsker's theorem we can express the limit process for $\sqrt{n%
}\left( \hat{F}_{y|x}-F_{y|x},\psi \right) $ as $(Q_{y|x},\psi )=Q_{\psi
}\left( U_{x},U_{xy}\right) $ by substituting the limit Browning bridge
processes for the arguments of $Q_{\psi }\left( .,.\right) .$

For any $\psi _{1},...,\psi _{l}\in D\left( W\right) $ the joint limit
process for 
\begin{equation*}
\sqrt{n}\left( \hat{F}_{y|x}-F_{y|x},\psi _{1}\right) ,...,\sqrt{n}\left( 
\hat{F}_{y|x}-F_{y|x},\psi _{l}\right)
\end{equation*}
is similarly given by the joint process of $Q_{\psi _{1}}\left(
U_{x},U_{xy}\right) ,...,Q_{\psi _{l}}(U_{x},U_{xy}).$ This is a Gaussian
process. The mean is zero since $Q_{\psi }$ is linear in its arguments and
the covariance is given by $cov\left( Q_{\psi _{1}}\left(
U_{x},U_{xy}\right) ,Q_{\psi _{2}}(U_{x},U_{xy})\right) =cov\left(
(Q_{y|x},\psi _{1}),(Q_{y|x},\psi _{2})\right) .$ Existence follows from
boundedness of the functions in the expressions and bounded support of $\psi
.$

By assumption of the theorem $h^{2}=o(n^{-\frac{1}{2}}),$ thus the limit
process is fully described by $Q_{y|x}.$

$\blacksquare $

Proof of Theorem 4.

For $\left( \ref{meanfunc}\right) $ we obtain%
\begin{eqnarray*}
&&\left( -1\right) ^{d_{x}+1}\int \int \hat{F}_{x,y}\left( x,y\right)
\partial ^{d_{x}}\psi \left( \hat{F}_{x}(x)\right) [y\psi _{v}^{\prime
}(y)+\psi _{v}\left( y\right) ]d\left( \hat{F}_{x}(x)\right) dy \\
&=&\left( -1\right) ^{d_{x}+1}\{\int \int F_{x,y}\left( x,y\right) \partial
^{d_{x}}\psi \left( F_{x}(x)\right) [y\psi _{v}^{\prime }(y)+\psi _{v}\left(
y\right) ]d\left( F_{x}(x)\right) dy \\
&&+\int \int [\hat{F}_{x,y}\left( x,y\right) -F_{x,y}(x,y)][y\psi
_{v}^{\prime }(y)+\psi _{v}\left( y\right) ]\partial ^{d_{x}}\psi \left(
F_{x}(x)\right) d\left( F_{x}(x)\right) dy \\
&&+\int \int F_{x,y}\left( x,y\right) \left( \partial ^{d_{x}}\right)
^{2}\psi \left( F_{x}(x)\right) [\hat{F}_{x}\left( x\right) -F_{x}(x)][y\psi
_{v}^{\prime }(y)+\psi _{v}\left( y\right) ]d\left( F_{x}(x)\right) dy \\
&&+\int \int F_{x,y}\left( x,y\right) \partial ^{d_{x}}\psi \left(
F_{x}(x)\right) [y\psi _{v}^{\prime }(y)+\psi _{v}\left( y\right) ]d\left( 
\hat{F}_{x}(x)-F_{x}\left( x\right) \right) dy \\
&&+\tilde{R}\},
\end{eqnarray*}%
where $\tilde{R}$ combines the remaining terms. Analogously to the proof of
Theorem 3 $\sqrt{n}\left( \hat{m}-m,\psi \psi _{v}\right) $ is represented
as 
\begin{equation*}
Q_{\psi \psi _{v}}\left( \sqrt{n}\left( \hat{F}_{x}-F_{x}\right) ,\sqrt{n}%
\left( \hat{F}_{xy}-F_{xy}\right) \right) +n^{-\frac{1}{2}}R\left( \sqrt{n}(%
\hat{F}_{x}-F_{x}),\sqrt{n}\left( \hat{F}_{xy}-F_{xy}\right) \right) .
\end{equation*}%
The limit process for the first functional is expressed via a value of the
functional for Brownian bridges, 
\begin{eqnarray}
Q_{\psi \psi _{v}}\left( U_{x},U_{xy}\right) &=&\int \int U_{x,y}[y\psi
_{v}^{\prime }(y)+\psi _{v}\left( y\right) ]\partial ^{d_{x}}\psi \left(
F_{x}(x)\right) d\left( F_{x}(x)\right) dy  \label{q} \\
&&+\int \int F_{x,y}\left( x,y\right) \left( \partial ^{d_{x}}\right)
^{2}\psi \left( F_{x}(x)\right) U_{x}[y\psi _{v}^{\prime }(y)+\psi
_{v}\left( y\right) ]d\left( F_{x}(x)\right) dy  \notag \\
&&+\int \int F_{x,y}\left( x,y\right) \partial ^{d_{x}}\psi \left(
F_{x}(x)\right) [y\psi _{v}^{\prime }(y)+\psi _{v}\left( y\right) ]d\left(
U_{x}\right) dy.  \notag
\end{eqnarray}

This process is Gaussian with mean zero; summing over $v$ we get a zero mean
limit process, $(Q_{m},\psi )=\Sigma _{v}Q_{\psi \psi _{v}}\left(
U_{x},U_{xy}\right) .$ We need to verify that the bilinear covariance
functional $cov\left( (Q_{m},\psi _{1}),(Q_{m},\psi _{2})\right) $ is
well-defined (bounded) for any $\psi _{1},\psi _{2}.$

Since expectation of $Q_{m}$ is zero%
\begin{eqnarray*}
\left\vert cov\left( (Q_{m},\psi _{1}),(Q_{m},\psi _{2})\right) \right\vert
&\leq &\left[ E(Q_{m},\psi _{1})^{2}E(Q_{m},\psi _{2})^{2}\right] ^{\frac{1}{%
2}}, \\
E(Q_{m},\psi )^{2} &=&E\left( \Sigma _{v}Q_{\psi \psi _{v}}\left(
U_{x},U_{xy}\right) \right) ^{2}.
\end{eqnarray*}%
Thus it is sufficient to consider variances for some $\psi .$

The representation in $\left( \ref{q}\right) $ involves three terms, it is
sufficient to show that the variance of the sum of each type of term over
all $v$ is bounded.

Recall that here $cov(U_{z_{1}},U_{z_{2}})=F(\tilde{z})-F\left( z_{1}\right)
F(z_{2}),$ where $\tilde{z}=z_{1}\wedge z_{2}.$

Start with the first term in $\left( \ref{q}\right) $ and consider its
variance$.$

Evaluate

\begin{eqnarray*}
&&E\{\int ...\int U_{x_{1},y_{1}}U_{x_{2},y_{2}}[y_{1}\psi _{v_{1}}^{\prime
}(y_{1})+\psi _{v_{1}}\left( y_{1}\right) ][y_{2}\psi _{v_{2}}^{\prime
}(y_{2})+\psi _{v_{2}}\left( y_{2}\right) ]dy_{1}dy_{2} \\
&&\cdot \partial ^{d_{x}}\psi \left( F_{x}(x_{1})\right) d\left(
F_{x}(x_{1})\right) \partial ^{d_{x}}\psi \left( F_{x}(x_{2})\right) d\left(
F_{x}(x_{2})\right) \} \\
&=&E_{1}-E_{1,2}\text{ with} \\
E_{1} &=&\{\int ...\int F(x_{1},y_{1})[y_{1}\psi _{v_{1}}^{\prime
}(y_{1})+\psi _{v_{1}}\left( y_{1}\right) ]\left[ \int^{y_{1}}[y_{2}\psi
_{v_{2}}^{\prime }(y_{2})+\psi _{v_{2}}\left( y_{2}\right) ]dy_{2}\right]
dy_{1} \\
&&\cdot \partial ^{d_{x}}\psi \left( F_{x}(x_{1})\right) d\left(
F_{x}(x_{1})\right) \int^{x_{1}}\partial ^{d_{x}}\psi \left(
F_{x}(x_{2})\right) d\left( F_{x}(x_{2})\right) \} \\
\text{and }E_{1,2} &=&\tilde{E}_{1}\tilde{E}_{2}\text{ where for }i=1,2 \\
\tilde{E}_{i} &=&\int ...\int F(x_{i},y_{i})[y_{i}\psi _{v}^{\prime
}(y_{i})+\psi _{v}\left( y_{i}\right) ]\partial ^{d_{x}}\psi \left(
F_{x}(x_{i})\right) d\left( F_{x}(x_{i})\right) dy_{i}.
\end{eqnarray*}

For $E_{1}$ integrating we get (dropping the subsript 1 on variables)%
\begin{equation*}
\int ...\int F(x,y)[y^{2}\psi _{v_{1}}^{\prime }(y)\psi _{v_{2}}(y)+y\psi
_{v_{1}}(y)\psi _{v_{2}}(y)]dy\cdot \frac{1}{2}\partial ^{d_{x}}\psi
^{2}\left( F_{x}(x)\right) dF(x).
\end{equation*}%
By construction of the partition of unity $\left\vert \Sigma \psi
_{v_{1}}^{\prime }(y)\psi _{v_{2}}(y)\right\vert $ as well as $\Sigma \psi
_{v_{1}}(y)\psi _{v_{2}}(y)$ are uniformly bounded, say both by some $\tilde{%
B}.$ We get 
\begin{eqnarray*}
&&\left\vert \int ...\int F(x,y)[y^{2}\psi _{v_{1}}^{\prime }(y)\psi
_{v_{2}}(y)+y\psi _{v_{1}}(y)\psi _{v_{2}}(y)]dy\cdot \frac{1}{2}\partial
^{d_{x}}\psi ^{2}\left( F_{x}(x)\right) dF(x)\right\vert \\
&\leq &\frac{\tilde{B}}{2}[\left( E_{|x}\left( y^{2}\right) ,\psi
^{2}\right) +\left\vert \left( E_{|x}y,\psi ^{2}\right) \right\vert ].
\end{eqnarray*}

Note that $\psi ^{2}\in D(W).$ By Assumption 5 then this contribution to the
covariance is bounded.

Similarly boundedness of the othe contributions from all the terms into the
covariance can be obtained. By the condition $h^{2}=o(n^{-\frac{1}{2}})$ on
the bandwidth the bias does not affect the limit process.

$\blacksquare $

\end{document}